\documentclass[12pt,reqno]{amsart}
\usepackage{times}
\usepackage{amssymb}
\usepackage{fullpage} 
\usepackage{setspace}
\usepackage{mathtools}
\usepackage{comment}
\usepackage{color}
\usepackage{amsmath,amsfonts,amsthm,amssymb,amscd,url}
\usepackage{graphicx}
\usepackage{mathrsfs}
\usepackage{todonotes}
\usepackage{dsfont}
\usepackage{enumitem}
\usepackage[parfill]{parskip} 
\usepackage[toc,page]{appendix}
\usepackage[utf8]{inputenc}
\usepackage{fullpage,amssymb,amsmath}
\usepackage{xcolor}
\usepackage{halloweenmath}
\usepackage{url}
\usepackage[colorlinks=true, pdfstartview=FitH, linkcolor=blue, citecolor=blue, urlcolor=blue]{hyperref}

\newcommand{\li}{\mathop{\rm li}}
\renewcommand{\mod}[1]{{\ifmmode\text{\rm\ (mod~$#1$)}\else\discretionary{}{}{\hbox{\!\!}}\rm(mod~$#1$)\fi}}

\newcommand{\fA}{{\mathfrak A}}
\newcommand{\bQ}{{\mathbb Q}}

\title{Comparative Prime Number Theory Problem List}
\author{Alia Hamieh, Habiba Kadiri, Greg Martin, and Nathan Ng}
\date{Last updated \today}

\begin{document}

\maketitle

\section*{Introduction}
This list of open problems related to the field of comparative prime number theory was created by the Pacific Institute of Mathematical Sciences (PIMS) Collaborative Research Group
\href{https://www.pims.math.ca/collaborative-research-groups/lfunctions}{$L$-functions in Analytic Number Theory} (2022--2025). 
The main focuses of this CRG are moments of $L$-functions and automorphic forms, explicit results in analytic number theory, and comparative prime number theory.
This open problem list was circulated in connection with the \href{https://sites.google.com/view/crgl-functions/comparative-prime-number-theory-symposium}{Comparative Prime Number Theory Symposium} (CPNTS) which took place at the University of British Columbia (Vancouver) from June 17--21, 2024.

For an overview of the field of comparative prime number theory, we recommend the introductory pages of {\em \href{https://arxiv.org/abs/2309.08729}{An annotated bibliography for comparative prime number theory}} by Martin {\em et~al.}~\cite{ABCPNT}. The field includes (but is not limited to) these topics:
\begin{itemize}\setlength\itemsep{-6pt}
\item classical prime-counting functions
\item races among primes in arithmetic progressions
\item races associated with elliptic curves, number fields, and function fields
\item summatory functions of arithmetic functions such as $\mu(n)$, $(-1)^{\Omega(n)}$, $\tau(n)$, etc.
\item the distribution of zeros of $L$-functions associated with the above races, including the linear independence of their imaginary parts
\item general oscillations and frequency of sign changes of number-theoretic error terms.
\end{itemize}

Our hope is that this list will stimulate research and lead to future collaborations among symposium participants and other researchers. In particular, problems on the list are available for anyone to decide to work on. If you make progress on any of the problems from this list, we would appreciate it if you acknowledged and cited this list in your work.


\section*{Notation}

As usual, $\pi(x)$ denotes the number of primes up to~$x$, and we have the weighted variants
\[
\theta(x) = \sum_{p\le x} \log p \qquad\text{and}\qquad \psi(x) = \sum_{n\le x} \Lambda(n) = \sum_{p^k\le x} \log p.
\]
The same functions have variants for primes in arithmetic progressions, such as
\[
\pi(x;q,a) = \#\{p\le x\colon p\equiv a\mod q\}.
\]
We write $\delta_{q;a_1,a_2,\ldots,a_r}$ for the logarithmic density of the set
\[
\bigl\{ x\ge2\colon \pi(x;q,a_1) > \pi(x;q,a_2) > \cdots \pi(x;q,a_r) \bigr\}.
\]
The Mertens sum is the summatory function $M(x) = \sum_{n\le x} \mu(n)$ of the Möbius function, and $M_0(x)$ is the modified Mertens sum where $\mu(x)$ is counted with weight~$\frac12$ if~$x$ is an integer.

We write $\chi_{-4}(n)$ for the nonprincipal Dirichlet character modulo~$4$. We also let $\omega(n)$ denote the number of distinct prime factors of~$n$ and $\Omega(n)$ the number of prime factors of~$n$ counted with multiplicity. We write $\gamma_1\approx 14.135$ for the imaginary part of the first nontrivial zero of $\zeta(s)$.

RH denotes the Riemann hypothesis, the conjecture that all nontrivial zeros of $\zeta(s)$ have real part equal to~$\frac12$. GRH denotes the generalized Riemann hypothesis, which is the same assertion for Dirichlet $L$-functions. SZ denotes the simple zeros conjecture, asserting that all zeros (of the $L$-function, including zeta-functions under discussion) are simple. LI denotes the linear independence conjecture, asserting that the nonnegative imaginary parts of the nontrivial zeros (of the $L$-function or family of $L$-functions under discussion) are linearly independent over the rational numbers.

\section*{Theoretical problems}

\begin{enumerate}[label=\#\arabic*.,leftmargin=0mm]

\item (Greg Martin and Lior Silberman, University of British Columbia)
An {\em ordinate} of an $L$-function is a positive imaginary part of a nontrivial zero of that $L$-function. LI implies that at most one positive ordinate can be rational, and presumably even that is not possible (for classical analytic $L$-functions).

{\bf Problem~1}: Prove that $\zeta(s)$ has at least one(!)~irrational ordinate.

{\bf Problem~2}: Prove that the positive ordinates of $\zeta(s)$ are not all rational multiples of one another.

These two problems can be extended to other $L$-functions (such as Dirichet $L$-functions) or even to collections of $L$-functions (such as the set of all Dirichlet $L$-functions corresponding to primitive characters; for example, ``prove that the positive ordinates of all Dirichlet $L$-functions are not all rational multiples of one another'').

One needs to be careful with higher-degree $L$-functions since there are known barriers to LI, such as: the factorization of Dedekind $\zeta$-functions and Artin $L$-functions into other $L$-functions; the vanishing of elliptic curve $L$-functions at the central point; and the ability to shift an entire $L$-function vertically if we don't restrict to the central character.

\item (Alexandre Bailleul, Université Paris--Saclay) As a concrete instance of the above problem:

{\bf Problem}: Show that $\gamma_1$ is irrational.

\item (Nicol Leong, University of New South Wales--Canberra at ADFA)

{\bf Problem}: Can we obtain upper bounds for $1/{|\zeta'(\rho)|}$ as $\rho$ runs over the nontrivial zeros of $\zeta(s)$?
 
Assuming SZ for $\zeta(s)$, the Mertens sum has the explicit formula
\begin{equation*}
M_0(x) = \sum_\rho \frac{x^\rho}{\rho \zeta'(\rho)} -2 + \sum_{n= 1}^\infty \frac{(-1)^{n-1}(2\pi)^{2n}}{(2n)!n\zeta(2n+1)x^{2n}}.
\end{equation*}
The order of the quantity $1/|\zeta'(\rho)|$ is a huge open problem that is difficult to approach. We have an exact formula (see~\cite{Tao})
\begin{equation}
  \label{taoformula}
\frac{1}{\zeta'(\rho)} = \lim_{s\to \rho}\frac{s-\rho}{\zeta(s)}= -\frac{2e^{1-B\rho}\rho(\rho-1)\Gamma(1+\rho/2)\pi^{-\rho/2}}{\prod_{\rho'\neq \rho}(1-\tfrac{\rho}{\rho'})e^{\rho/\rho'}},
\end{equation}
where $B= -0.0230957\ldots$ is the constant appearing in the Weierstrass product for $\xi(s)$. Thus the main issue here is bounding the term $1- \tfrac{\rho}{\rho'}$, in other words, obtaining a lower bound for the distance $|\rho -\rho'|$ from a zero to its nearest neighbour. Even assuming SZ, we cannot yet rule out the existence of two extremely close zeros.

\item (Nathan Ng, University of Lethbridge)

{\bf Problem}: Use the formula~\eqref{taoformula} to compute an upper bound for 
\begin{equation}
  \label{negmoment}
 \sum_{0 < \gamma < T} \frac{1}{|\zeta'(\rho)|^{2k}}
\end{equation}
for positive real numbers~$k$ (the cases $k=\frac12$ and $k=1$ are especially interesting), assuming SZ and the GUE conjecture on the nearest-neighbour spacing on zeros of $\zeta(s)$.

At the moment there are no upper bounds (even conditionally) for this exact expression.
The recent article of Bui, Florea, and Milinovich~\cite{BFM} gives some first nontrivial upper bounds, but for sums of the shape 
\[
    \sum_{\substack{0 < \gamma < T \\ \gamma \in S}} \frac{1}{|\zeta'(\rho)|^{2k}}
\]
where~$S$ is a particular subset of the zeros of $\zeta(s)$ on which $\zeta'(\rho)$ is not too small. It is believed that the behaviour of the sum~\eqref{negmoment} changes once $k\ge\frac32$; we know of a graduate student who is working on a suitable model to make conjectures in this region.

\item (Nathan Ng, University of Lethbridge) My Ph.D.~thesis \cite{NgPhd} contains the conjecture
\[
      \sum_{0 < \gamma < T}  \biggl| \frac{\zeta(2 \rho)}{\zeta'(\rho)} \biggr|^2
      \sim \frac{1}{2 \pi} T.
\]
{\bf Problem 1}:  Use the argument in Milinovich and Ng~\cite{MN} to establish 
the lower bound 
\[
    \sum_{0 < \gamma < T} \biggl| \frac{\zeta(2 \rho)}{\zeta'(\rho)} \biggr|^2 \ge \biggl( \frac{1}{4 \pi} + o(1) \biggr) T.
\]
{\bf Problem 2}:  Prove lower bounds for the discrete moments of the shape
\[
      \sum_{0 < \gamma < T} \biggl| \frac{\zeta(2 \rho)}{\zeta'(\rho)} \biggr|^{2k}
      \gg T (\log T)^{(k+1)^2}.
\]
See the articles of Heap, Li, and Zhao~\cite{HLZ} and Gao and Zhao~\cite{GZ}, where the authors prove lower bounds of the type
\[
      \sum_{0 < \gamma < T}  \frac{ 1}{|\zeta'(\rho)|^{2k}} 
      \gg T (\log T)^{(k+1)^2}
\]
for $k >0$.

\item (Nicol Leong, University of New South Wales--Canberra at ADFA)

{\bf Problem}: Can we improve upper bounds for the multiplicities of nontrivial zeros of $\zeta(s)$?
 
Let $m(\rho)$ denote the multiplicity of a nontrivial zero of $\zeta(s)$. We know unconditionally that
\begin{equation*}
m(\beta +i\gamma) \ll \log \gamma,
\end{equation*}
while under the assumption of the Lindel\"of hypothesis or the Riemann hypothesis, the bound can be improved to
\begin{equation*}
o(\log \gamma) \qquad \text{or} \qquad \ll\frac{\log \gamma}{\log\log \gamma},
\end{equation*}
respectively.
Ivi\'c has proved a couple of unconditional results of order $\ll \log\gamma$, but this is still far from Karatsuba's conjectured size of $\ll 1$. (Of course, under the LI conjecture, $m(\rho)$ would equal~$1$ for all~$\rho$.) For more details on the problem and Ivi\'c's approach, see for instance~\cite{Ivic1,Ivic2}.

\item (Alexandre Bailleul, Université Paris--Saclay) Skewes's number is defined as the smallest~$x$ such that $\pi(x) > \li(x)$. It was shown to exist by Littlewood, and a first (enormous) bound was found by Skewes. Its value is still unknown, but it must be quite large since it is known to lie in the interval $(10^{19}, 1.3971653 \times 10^{316})$ (the lower bound is due to Büthe~\cite{Buthe0}, while the upper bound is due to Saouter, Trudgian, and Demichel~\cite{STD}); experts suggest that it is quite close to this upper bound. The corresponding logarithmic density for the race between $\pi(x)$ and $\li(x)$ was determined to be extremely small by Rubinstein and Sarnak~\cite{RS}, namely $\delta(\pi,\li) \approx 2.6 \times 10^{-7}$.

Another example of a race with a large bias is the race between primes congruent to $2\mod3$ versus primes congruent to $1\mod3$. Rubinstein and Sarnak~\cite{RS} calculated the corresponding logarithmic density $\delta_{q;1,2} \approx 1.0\times10^{-3}$, while Bays and Hudson shows that the first sign change of $\pi(x;3,2) - \pi(x;3,1)$ occurs at $x=608{,}981{,}813{,}029$. Therefore, it seems natural to think that the size of the ``Skewes's number'' of a prime number race is correlated with how biased that race is. In particular, classical prime number races have their logarithmic densities tending to~$\frac12$ as $q \to \infty$, and so we might expect the corresponding Skewes's numbers to be less and less extreme (but of course bounded below by the modulus~$q$).

{\bf Problem}: Obtain an upper bound for the first sign change of $\pi(x;q,a)-\pi(x;q,b)$ in terms of~$q$, preferably one that increases as slowly as possible.

\item \label{johnston1} (Daniel R.~Johnston, University of New South Wales--Canberra at ADFA) This question was posed during Ethan Lee's talk at the CPNTS.

Given a number field $\mathbb{K}$ with ring of integers $\mathcal{O}_{\mathbb{K}}$, let $\kappa_{\mathbb{K}}$ be the residue at the pole $s=1$ for the Dedekind zeta-function $\zeta_{\mathbb{K}}(s)$ associated to $\mathbb{K}$. Consider the difference
\[
    \Delta_{\mathbb{K}}(x)=\prod_{N(\mathfrak{p})\leq x}\left(1-\frac{1}{N(\mathfrak{p})}\right)^{-1}-e^{\gamma}\kappa_{\mathbb{K}}\log x,
\]
where the product runs over the prime ideals $\mathfrak{p}$ of $\mathcal{O}_{\mathbb{K}}$ with norm $N(\mathfrak{p})$ up to~$x$; this difference $\Delta_K(x)$, which tends to be positive, is connected to Mertens's third theorem for number fields. B\"uthe~\cite{buthe2015first} showed that the first negative value of $\Delta_{\mathbb{Q}}(x)$ occurs before $2\cdot 10^{215}$. On the other hand, as shown in~\cite[Table~1]{hathi2021mertens}, the tendency for $\Delta_{\mathbb{K}}(x)$ to be positive is weaker (in terms of logarithmic density) for the simple number fields $\mathbb{K}=\mathbb{Q}(\sqrt{5})$ and $\mathbb{K}=\mathbb{Q}(\sqrt{13})$ than for $\mathbb{K}=\mathbb{Q}$. It's unclear whether this should mean that the first negative values of $\Delta_{\mathbb{K}}(x)$ in these cases should occur sooner than those of $\Delta_{\mathbb{Q}}(x)$.

{\bf Problem:} Localize the first negative value for $\Delta_{\mathbb{K}}(x)$ when $\mathbb{K}$ is $\mathbb{Q}(\sqrt{5})$, $\mathbb{Q}(\sqrt{13})$, or some other simple number field. In general, obtain an upper bound for the first sign change of $\Delta_{\mathbb{K}}(x)$ for a family of number fields~$\mathbb{K}$.

\item (Alexandre Bailleul, Université Paris--Saclay) There exist highly biased races, such as the (appropriately normalized) race between quadratic residues and nonresidues for composite moduli~$q$. Fiorilli~\cite{Fio} formulated a conjecture for the Skewes's number of this race between residues and nonresidues (the smallest real number for which the counting function for primes that are quadratic residues overcomes the bias): it is roughly of the form $\exp\bigl( \exp\bigl( {\rho(q)}/\log{\mathrm{rad}(q)} \bigr) \bigr)$, where $\rho(q)$ is the index of the subgroup of squares in $(\mathbb Z/q\mathbb Z)^{\times}$ and $\mathrm{rad}(q)$ is the radical (largest squarefree divisor) of~$q$.

More recently, Fiorilli and Jouve~\cite{FJ} and Bailleul~\cite{Bai} have investigated similar highly-biased races in number fields.

{\bf Problem}: Bound Skewes's numbers of other highly biased races in terms of the parameters defining the race, either conjecturally or provably.

\item (Alexandre Bailleul, Université Paris--Saclay) Consider Chebyshev's bias in the function field setting. As is customary since the work of Cha~\cite{Cha}, in addition to assuming the linear independence (over~$\bQ$) of the corresponding arguments of zeros of $L$-functions, one typically assumes that those arguments are linearly independent with~$\pi$ as well. While there are no known rational linear relations among the arguments of these zeros themselves, there are known instances where such an argument is a rational multiple of~$\pi$. These instances usually come from special geometric conditions, such as supersingularity of elliptic curves or jacobians that are not absolutely simple.

{\bf Problem}: Show that the failure of linear independence with $\pi$ always has a geometric explanation.

\item (Stanley Yao Xiao, University of Northern British Columbia)
Let $f \in \mathbb{Z}[x,y]$ be a primitive, positive definite binary quadratic form such that $f(x,1) \not \equiv x(x+1) \pmod{2}$. Consider the counting function for the generalized Friedlander--Iwaniec primes
\[
\pi_f(x) = \# \{p \leq x \colon p  = f(a,b^2) \text{ for some integers $a$ and $b$}\}.
\]
By the main theorem in my paper~\cite{X}, we have $\pi_f(x) \sim C(\operatorname{Disc}(f)) x^{3/4}/\log x$ for some positive constant~$C(\operatorname{Disc}(f))$ depending only on the discriminant of~$f$. 

{\bf Problem}: If~$f$ and~$g$ are two distinct such quadratic forms with the same discriminant, how do the quantities $\pi_f(x)$ and $\pi_g(x)$ compare as $x \rightarrow \infty$?

A similar question can be asked for indefinite binary quadratic forms, where the analogous counting function must restrict~$a$ and~$b$ to a suitable ``box'' depending on~$x$.

Heath-Brown and Moroz~\cite{HBM} gave a similar asymptotic formula for primitive binary cubic forms~$F$ such that $F(x,y) \not \equiv xy(x+y) \pmod{2}$; one could again look for pairs $F,G$ of distinct such cubic forms such that $\pi_F(x) \sim \pi_G(x)$ and compare those counting functions to each other.

\item (Shin-ya Koyama, Toyo University)
A new formulation of Chebyshev's bias was given in~\cite{AK} in terms of the difference of the weighted counting function of prime numbers~$p$ as
\begin{equation*}
\sum_{\genfrac{}{}{0pt}{1}{p\le x}{p\equiv 3\mod4}} \frac 1{\sqrt p}
-\sum_{\genfrac{}{}{0pt}{1}{p\le x}{p\equiv 1\mod4}} \frac 1{\sqrt p}
\sim\frac12\log\log x\quad(x\to\infty).
\end{equation*}
{\bf Problem}: What is the logical inclusion or relation of this new definition to the conventional definiton of ``Chebyshev's bias'' in terms of the logarithmic density of the set of~$x$ satisfying $\pi(x,4,3)>\pi(x,4,1)$ being very close to 1?

\item (Greg Martin, University of British Columbia)
This question was asked by Carl Pomerance and appears in~\cite[Section~3]{Mar}. Let~$a$ and~$b$ be distinct reduced residue classes modulo~$q$.

{\bf Problem}: Prove unconditionally that the set of~$x$ for which $\pi(x;q,a)=\pi(x;q,b)$ (the set of ``ties'' in the prime number race) has (natural or logarithmic) density~$0$.

If we assume GRH and a very mild version of LI (see~\cite[Theorem~2.2]{Devin}), the limiting logarithmic distribution of $\bigl( \pi(x;q,a)=\pi(x;q,b) \bigr) (\log x)/\sqrt x$ exists and is absolutely continuous with respect to Lebesgue measure. Consequently, for any function $f(x)$ that is $o(\sqrt x/\log x)$, the set of~$x$ for which $\bigl| \pi(x;q,a)=\pi(x;q,b) \bigr| < f(x)$ has logarithmic density~$0$; this conclusion is far stronger than ``ties have density~$0$''. But ``ties have density~$0$'' is such a weak assertion that we might hope to prove it unconditionally.

\end{enumerate}

\section*{Explicit versions of theoretical results}

\begin{enumerate}[resume,label=\#\arabic*.,leftmargin=0mm]

\item (Daniel R.~Johnston, University of New South Wales--Canberra at ADFA)
Define
\[
\Delta^{\theta_r}(x) = \sum_{p\leq x}\frac{\log p}{p}-(\log x+E),
\quad\text{where }
E=-\gamma-\sum_p\frac{\log p}{p(p-1)}=-1.33258\ldots,
\]
to be the error term in Mertens's first theorem. Computations~\cite[Theorem~21]{rosser1962approximate} have shown that $\Delta^{\theta_r}(x) > 0$ for all $0<x\leq 10^8$, but Lay~ \cite{lay2015sign} proved that $\Delta^{\theta_r}(x)$ changes sign infinitely often.

{\bf Problem}: Find an explicit constant~$X$ such that $\Delta^{\theta_r}(x)<0$ for some $x\le X$.

An analogous upper bound has already been found~\cite{buthe2015first} for the first sign change in Mertens's second theorem, which implies the same upper bound for the first sign change in Mertens's third theorem as mentioned in problem~\ref{johnston1}

\item (Chiara Bellotti, University of New South Wales--Canberra at ADFA)

Various authors \cite{DvS,Hur,Ram,Sch} have given explicit upper bounds for the Mertens sum of the form $Ax/(\log x)^B$.
The current best known non-explicit unconditional bound for the Mertens sum (see \cite[pp.~309--315]{Ivic3}) is of the form
\[
M(x) \ll x \exp \bigl(-C (\log x)^{{3}/{5}} (\log \log x)^{-{1}/{5}}\bigr).
\]
{\bf Problem}: Find an explicit version of this upper bound, which would thus give a sharper bound for $M(x)$ for large~$x$.

\item (Daniel R.~Johnston, University of New South Wales--Canberra at ADFA)
Define
\[
S(x)=\sum_{n\leq x}(-2)^{\Omega(n)}
\qquad\text{and}\qquad
\alpha = \limsup_{x\to\infty}\frac{|S(x)|}{x}.
\]
In a soon-to-be-released preprint by students at UNSW--Canberra, it is shown that~$\alpha$ exists, and a (large) upper bound is given. Note that $S(x)$ has jump discontinuities of size~$x$ when~$x$ is a power of~$2$, which implies the lower bound $\alpha\ge\tfrac{1}{2}$.
A conjecture \cite[Conjecture~1.1]{sun2016pair}, supported by extensive computations \cite{mossinghoff2021oscillations}, asks whether $\alpha\leq 1$.

{\bf Problem}: What is the exact value of~$\alpha$?

\item (Greg Martin, University of British Columbia)
The standard way to investigate the frequency of sign changes of error terms such as $\psi(x)-x$ is to use repeated averaging so that the largest term in the explicit formula eventually dominates the sum of all other terms. In this way, one shows that the number of sign changes is essentially at least as large as the number of sign changes of that largest term. For example, Kaczorowski~\cite{K1} showed unconditionally that
\[
\liminf \frac{W^\psi(T)}{\log T} \ge \frac{\gamma_1}{\pi}.
\]
This result can be generalized to many other error terms (although more generally one must assme that the relevant $L$-functions have no real zeros, and GRH tends to be necessary for error terms with biases).

The best improvement to this lower bound I am aware of is by Morrill, Platt, and Trudgian~\cite{MPT} (improving work of Kaczorowski~\cite{K2}), who proved that
\[
\liminf \frac{W^\psi(T)}{\log T} \ge \frac{\gamma_1}{\pi} + 1.867\cdot10^{-30}.
\]

{\bf Problem 1}: Improve this lower bound, even assuming RH. Here are some suggested target lower bounds for $\liminf {W^\psi(T)}/{\log T}$:
(a) $1.01 \frac{\gamma_1}\pi$;
(b) $2 \frac{\gamma_1}\pi$;
(c) $\infty$.

What we expect to be true is that the large-scale intervals of positive and negative values of $\psi(x)-x$ do indeed alternate at about the speed given by the repeated-averaging argument, but that each such sign change detected is actually a flurry of many sign changes on a small scale (this is how random walks behave, for example). In particular, we expect $W^\psi(T)$ and related sign-change counting functions to grow more like $\sqrt T$ (up to logarithmic factors) than like $\log T$.

{\bf Problem 2}: Make any progress whatsoever towards advancing our knowledge of sign changes beyond the repeated-averaging techniques.

\item (Alexandre Bailleul, Université Paris--Saclay) This question was raised by Nilotpal Sinha~\cite{Sinha}. Let $P^+(n)$ denote the largest prime factor of~$n$. Numerical investigation suggests that there is a bias towards integers~$n$ for which $P^+(n)\equiv3\mod4$ over integers~$n$ for which $P^+(n)\equiv1\mod4$, but the asymptotic behaviour of $\sum_{n\le x} \chi_{-4}(P^+(n))$ is still unknown.

{\bf Problem}: Prove or disprove that $\sum_{n\le x} \chi_{-4}(P^+(n)) < 0$ for all $n \geq 2$. Possibly even find an asymptotic formula for the sum.

\end{enumerate}

\section*{Computational problems}

\begin{enumerate}[resume,label=\#\arabic*.,leftmargin=0mm]

\item (Nathan Ng, University of Lethbridge)
Feuerverger and Martin~\cite{FeM} numerically computed the densities for some three-way races, as well as for the four-way races for the moduli~$5$, $8$, and~$12$. Their paper contains a general formula for densities in races with any number of competitors; but there have been few (if any) subsequent computations of densities for races with three or more competitors.

{\bf Problem 1}: Numerically compute some further three-way (or higher) races.

Fiorilli and Martin~\cite{FiM} determined the 117 most biased two-way races.

{\bf Problem 2}: Determine the $100$ (for example) most biased three-way races.

It could be more challenging to bound the densities in three-way races sufficiently well to verify that such a list is exhaustive: as Lamzouri~\cite{Lam} showed, the rate of convergence of three-way race densities to $\frac16$ is only $O(1/\log q)$, as opposed to a $q^{-1/2+o(1)}$ rate of convergence to $\frac12$ for densities in the two-way case.

\item (Greg Martin, University of British Columbia)

{\bf Problem}: Compute (assuming GRH and LI) the densities in the five-way race among the quadratic residues and among the quadratic nonresidues modulo~$11$.

There are $120$ densities of the form $\delta_{11;a_1,a_2,a_3,a_4,a_5}$ where $\{a_1,a_2,a_3,a_4,a_5\} = \{1,3,4,5,9\}$ (the quadratic residues modulo~$11$). There is a known formula~\cite[Theorem~4]{FeM} for $\delta_{q;a_1,a_2,a_3,a_4,a_5}$ that can in principle be computed to any desired accuracy; however, the formula involves a four-dimensional integral that was beyond the authors' computational ability at the time.

In addition to being the simplest natural five-way prime number race, Bays and Hudson~\cite{BH} noticed that this particular race exhibited an interesting cyclic phenonemon, which the computation of these densities might shed some light on.

Because of the symmetries in~\cite[Theorem~2]{FeM}, there are only $8$ (presumably) distinct values among these $120$ densities, so all that would need to be computed are $\delta_{11;1, 3, 4, 5, 9}$, $\delta_{11;1, 3, 4, 9, 5}$, $\delta_{11;1, 3, 5, 4, 9}$, $\delta_{11;1, 3, 5, 9, 4}$, $\delta_{11;1, 3, 9, 5, 4}$, $\delta_{11;1, 5, 3, 4, 9}$, $\delta_{11;1, 5, 3, 9, 4}$, and $\delta_{11;1, 5, 9, 3, 4}$. By the same theorem, these values also represent all possible densities in the five-way race among the quadratic nonresidues $\{2,6,7,8,10\}$ modulo~$11$.

\item (Nathan Ng, University of Lethbridge)
My Ph.D.~thesis \cite{NgPhd} contains some computations of numerical examples of prime number races in Chebotarev's density theorem. There is now a large database of zeros of Artin $L$-functions at the LMFDB~\cite{LMFDB}, as well as Rubinstein's lcalc program which has been integrated into SageMath.

{\bf Problem}: Investigate further examples. Determine the most biased two-way Chebotarev races.

\item (Greg Martin, University of British Columbia)
In 1959 Shanks~\cite{Shanks} made an observation that seems to have been overlooked since then: some combinatorial reasoning led him to assert that for each integer $k\ge1$, the odd numbers~$n$ with $\Omega(n)=k$ should be biased towards the residue class $(-1)^k\mod4$, where the case $k=1$ is Chebyshev's original observation. In other words, $(-1)^{\Omega(n)}\chi_{-4}(n)$ should be biased towards~$1$ over~$-1$.

{\bf Problem}: Assuming GRH and LI, calculate the logarithmic density of the sets
\[
\{n\ge1\colon (-1)^{\Omega(n)}\chi_{-4}(n) = 1\} \quad\text{and}\quad \{n\ge1\colon (-1)^{\Omega(n)}\chi_{-4}(n) = -1\}.
\]
Presumably a similar phenomenon should hold for $(-1)^{\Omega(n)}\chi(n)$ for any quadratic character~$\chi$; the case $\chi=\chi_{-3}$ (the nonprincipal Dirichlet character modulo~$3$) might be particularly pronounced.

By the way,~\cite{Shanks} also contains the timeless quip ``Since $\log\log\log n\to\infty$ (with great dignity)\dots.''

\item (Daniel R.~Johnston, University of New South Wales--Canberra at ADFA).
It has been shown \cite[Theorem~1.1]{johnston2023average} that RH is equivalent to
\[
\fA_1^\pi(x) = \int_0^x \bigl( \pi(t)-\li(t) \bigr) \,dt < 0 \quad\text{for all $x>2$.}
\]
Roughly, the reason this inequality holds under RH is that the constant $\sum_{\rho} 1/\rho(\rho+1)$ occuring in the explicit formula for $\fA_1^\pi(x)$ is too small to overcome the bias. However, Martin pointed out (private communication) that when looking at the analogous
\[
\fA_1^\pi(x;q,a) = \int_0^x \bigl( \phi(q)\pi(t;q,a)-\li(t) \bigr) \,dt,
\]
the corresponding sum over zeros of Dirichlet $L$-functions grows faster than the bias as~$q$ increases, and thus at some point the inequality $\fA_1^\pi(x;q,a) < 0$ should fail when~$q$ is large.

{\bf Problem}: Assuming GRH, classify those values of~$q$ and~$a$ such that $\fA_1^\pi(x;q,a)<0$ for all sufficiently large~$x$.

\end{enumerate}



\end{document}